\documentclass[11pt]{article}
\usepackage{indentfirst,latexsym}

\renewcommand{\normalsize}{\fontsize{10}{12}\selectfont}
\newcommand{\+}{\fontsize{11}{12}\selectfont}

\DeclareMathSizes{10}{12}{11}{10}
\DeclareMathSizes{9}{10}{8}{7}
\DeclareMathSizes{8}{10}{7}{7}
\DeclareMathSizes{14.4}{14.4}{11}{10}
\DeclareMathSizes{12}{14.4}{11}{10}

\newcommand{\<}{\hspace{-0.5em}}

\newcommand{\pdf}[0]{{\it pdf}\hspace{0.08 em}}

\newcommand{\D}{\mathrm{d}}

\begin{document}
\title {
\vspace{-1.3 cm}
Direct pivotal predictive inference
\thanks {{\normalsize \hspace{0.1 em}
With a typo correction and minor improvements
over earlier drafts (privately communicated).
\newline \indent {\small $^{\mathsection}$} \hspace{0.1 em}
{\tt \fontsize{11}{12}\selectfont anakreon@hol.gr}
\newline \indent {\small $^{\mathparagraph}$} \hspace{0.1 em}
2215 Coover Hall,
Iowa State University,
Ames, Iowa 50011;
{\tt \fontsize{11}{12}\selectfont berleant@iastate.edu}
}}
\vspace{-0.45 cm}
}
\author {
George Kahrimanis {\fontsize{8}{12}\selectfont $^{^{\mathsection}}$}
and
Daniel Berleant {\fontsize{8}{12}\selectfont $^{^{\mathparagraph}}$}
}
\date   {
\vspace{-0.22cm}
{\normalsize      19 October 2004   }}
\maketitle

\vspace{-0.9 cm}

\begin{abstract}
{\normalsize
Without assuming any \pdf\ for some measured parameter,
we derive a predictive \pdf\ for the outcome of a second measurement,
given the outcome of the first measurement
and two common assumptions about the noise.
These are that (1) it is additive, and (2) it is of some known \pdf.
The argument is based on a Bayesian analysis of the noise
when no \pdf\ is provided for the value of the parameter.
In this way we avoid assuming an {\it ad-hoc }prior.
We clarify how this method of direct predictive inference
is distinct from fiducial prediction.
We specify the distinct flaw in the fiducial argument,
and outline the importance of this development
in the foundations of probability and statistics.

\vspace{-0.2 cm}
}
\end{abstract}

\section{Outline and Motivation}
We examine the problem of deriving a \pdf\ 
predicting a future outcome $x_2$
based on an earlier outcome $x_1$
and a model of location measurements,\footnote
{{\normalsize \hspace{0.1 em}
`Location' means that the noise is additive
and its \pdf\ does not depend on the measured magnitude.
}}
when there is no prior \pdf.
We consider the solution based on the distribution of $x_2 {-}$ \< $x_1$,
which is readily determined.
In this way we avoid assuming some ``reference prior''
({\it e.g.} \cite{Kass_Wass}).
We also dodge the problematic premises of fiducial prediction
(\cite{Fisher35}, \cite{Fisher56}; reviewed in \cite {Seidenfeld92});
that is, we can do without parametric inference,
unlike in the Bayesian and fiducial treatments.

The intuitive appeal of the ``direct'' pivotal argument
(the term is borrowed from \cite{Seidenfeld92}, p.~365)
is balanced by an understandable skepticism toward pivotal inference
because an argument of this kind has also been advanced by Fisher
in support of {\em fiducial} inference ({\it e.g.} \cite{Fisher56}).
This generated much controversy in its day
and is generally regarded as problematic at best
(\cite{Pedersen}, \cite{Seidenfeld92}, \cite{Zabell}, \cite{Efron} ---
even by Fisher himself, in private communication
with Barnard and Savage, cited in \cite{Zabell}, p.~381).

In this setting we have a double task.
In the first place
we have to account for selectively dismissing fiducial inference
but not the pivotal argument of direct predictive inference.
We focus on a technical shortcoming {\em specific} to the fiducial argument,
in Sec.~\ref{sec:fiducial}.
This observation is interesting in itself,
because it complements earlier critical reviews of the fiducial argument,
like \cite{Pedersen}, \cite{Seidenfeld92}, and \cite{Zabell}.
On the other hand, to address certain well-founded concerns,
we submit a mathematical demonstration (starting at Sec.~\ref{sec:proof})
that the pivotal \pdf\ in the direct pivotal argument
is unaffected by our knowing $x_1$.

Doing so, we justify solutions to common practical problems,
without resorting to ``statistical principles'' and makeshift priors.
This treatment also introduces an important theoretical development,
because until now the fiducial argument has posed a profound quandary:
it has already been shown that it leads to contradiction
(see \cite{Seidenfeld92}, Sec.~5)
but, inasmuch as no specific step in the fiducial argument
has been identified as being wrong,
the foundation of probability theory (not only statistics)
is implicated, by default.
In the aftermath,
a tacit stopgap injuction has been generally in effect,
to avoid fiducial inference,
as it were with the force of an {\it ad-hoc }modification
of the axioms of probability.
In consequence, the process of deduction has been undermined,
because new axioms have not been explicitly stated,
and, even if they had, any modified foundation would seem contrived
in comparison to the classical one.
(This is most evident in the case of the ``direct'' pivotal argument,
which is similar to fiducial inference yet is distinct from it,
so that it is unclear whether the informal injunction applies.)
The present treatment dissolves this predicament,
by isolating the distinct flaw in the fiducial argument,
and so restoring the classical foundation in the theory of probability.

\section{The main argument}
\subsection{Fundamental issues related to predictive inference}
\label{sec:fundamental}
To use a plain example,
let $\theta$ be an unknown real constant,
not the result of any known random process.
Two independent measurements of $\theta$,
with outcomes $x_1$ and $x_2$,
correspond to conditional random variables
$X_{1\,\{\theta=t\}}$ and $X_{2\,\{\theta=t\}}$.
We assume that, conditionally on any possible value of $\theta$,
these are independent and distributed normally:
$X_{1\,\{\theta{=}t\}}$ $\sim$ $\mathrm{N}(t,\sigma_1^{\;2})$ and
$X_{2\, \{\theta{=}t\}}$ $\sim$ $\mathrm{N}(t,\sigma_2^{\;2})$,
where $\sigma_1$ and $\sigma_2$ are assumed constant and known.
The difference of the two outcomes, $d$ ${\equiv}$ $x_2{-}x_1$,
corresponds to $X_{2\,\{\theta{=}t\}} {-}$ \< $X_{1\,\{\theta{=}t\}}$,
which is distributed normally with mean $0$
and variance $\sigma_1^{\;2}{+}$ \< $\sigma_2^{\;2}$,
without reference to the value of $\theta$.
Therefore random variable $D$ can be defined,
following $\mathrm{N}(0,$ \< $\sigma_1^{\;2}{+}$ \<
$\sigma_2^{\;2})$,
and is not conditional on $\theta$.
Any random variable with this property
is called {\em pivotal} or a {\em pivot}.

The \pdf\ of $D$ can be directly applied to generate
the \pdf\ for ``$x_2$, given $x_1$''
according to what has been called
``direct'' pivotal argument (as distinct from fiducial prediction)
in \cite{Seidenfeld92} (p.~365).
But this is not the end of the issue,
because this argument relies on the premise
that the \pdf\ of $D_{\{x_1{=}s\}}$ is the same as the \pdf\ of $D$.
To defend this claim, it is not enough to know
that the distribution of $D$ is independent of $\theta$;
one must also demonstrate that knowledge of $x_1$
does not make any difference
(that is, $x_1$ by itself
does not specify any recognizable subset
of the reference class (or ``reference set'') associated with $D$
({\it e.g.} \cite{Fisher56} pp.~57-58, \cite{Pedersen},
or \cite{Seidenfeld92})).

The problem has persisted for decades, because,
lacking a mathematical demonstration of this claim
(\cite{Pedersen} Sec.~4; \cite{Zabell} Sec.~7.2),
conjectures about the post-data \pdf\ of a pivot
have been presumed by some authors
on the basis of an intuitive conviction,
starting with Fisher (e.g., see \cite{Zabell} Sec.~7.2),
who has evoked a general version of this assertion
in support of his fiducial argument.
The inadequacy of a bare appeal to intuition has been outstanding
since Fisher's assertion was checked wrong in a particular situation
\cite{Buehler_Fed} which involves the $t$-distribution 
(for an outline see \cite{Seidenfeld92} p.~364).

Nevertheless this finding does not extend to our example.
Moreover, in this work we submit a proof backing the claim
that, in the case of location measurements,
$x_1$ is irrelevant (by itself) to the \pdf\ of $D$.

\subsection{General considerations}
\subsubsection{Additive noise as a pivotal random variable}
Each location measurement (here we label them by $i = 1,2,...$)
can be thought of as involving a random process,
generating noise $E_i$ of a known \pdf\ $f_i(e_i)$,
which is then added to the unknown parameter $\theta$
to provide an outcome $x_i$.
Our assignment of a pre-data \pdf\ to $E_i$ is based on this assumption,
which we call Assumption~${\mathcal B}_i$,
or simply ${\mathcal B}_i$.

We shall employ a detailed notation for probability statements,
to display the defining assumptions (the ``context'' or "reference class")
which label the corresponding probability space.
For instance, to state explicitly that the pre-data \pdf\ of $E_1$
does not depend on (hypothetical values of) $\theta$, we write
\begin{equation} \label{eq:E_ind_Theta}
(\forall \, t) \ \ \ 
{\Pr}_{E_1}(e_1|
{\mathcal B}_1,\textrm{``}\theta{=}t\textrm{''})
\ \equiv \ 
{\Pr}_{E_1}(e_1|{\mathcal B}_1)
\ \equiv \ f_1(e_1) \, \textrm{.}
\end{equation}
(That is, $E_1$ is a pivot.)

The suffix in a probability statement has a double role.
Not only it denotes the random variable (here: $E_1$)
which is associated with the indicated value (here: $e_1$),
but also, when the random variable is continuous,
it specifies the parametrisation used to represent
the probability {\em density} function as a regular function.

In the following we shall focus on determining
the post-data \pdf\ of the error of measurement.
This \pdf\ can be expressed as ${\Pr}_{E_1}(e_1|
{\mathcal B}_1,\textrm{``}x_1{=}s\textrm{''})$.

\subsubsection{Disregarding the true value of the parameter}
\label{sec:disregarding}
In the expression ${\Pr}_{E_1}(e_1|
{\mathcal B}_1,\textrm{``}x_1{=}s\textrm{''})$,
the reference class is specified by two conditions:
an assumption regarding the noise of this type of measurement,
and an acceptance of the outcome of this particular measurement.

As a technical excercise,
if in the specification of the reference class we also included
the hypothetical ``true value of $\theta$'',
as in ${\Pr}_{E_1}(e_1|{\mathcal B}_1,
\textrm{``}x_1{=}s\textrm{''},\textrm{``}\theta{=}t\textrm{''})$,
the \pdf\ would collapse into the delta-function $\delta(e_1{-}s{+}t)$,
which is unspecified, therefore useless for predictive inference
(besides, it does not contain any trace
of the known properties of the measurement).

However, the general issue of selecting the reference class is still open;
at any rate,
it cannot be decided within the theory of probability.
We cannot {\em prove} that it is wrong to include ``$\theta{=}t$''
in the specification of the reference class.
We rather point out some counterintuitive consequences of this choice,
ultimately related to practical disadvantages.

A simple version of the same problem will arise if one tosses a fair coin,
and immediately covers it with a bowl.
Do we accept that the reference class consists of all such trials
(disregarding that either `heads' or `tails' has already become
a constituent of reality)
so that probability of heads-under-the-bowl is $0.5$,
or do we restrict the refererence class to this single case,
so that probability cannot be defined (except trivially)?
With the first option, we stand to gain (in the long run)
from betting against someone who wrongly believes that
probability of heads-under-the-bowl is $0.6$.
Not so if we follow the second option.

The same problem arises in any application of Bayes' theorem,
even if it is based on well-defined prior probability.
Take for instance the interpretation of a medical diagnostic test,
such as an HIV-antibody test.
The lab procedure outputs the relative likelihood of infection,
which then can be combined with prior probability based on information
about the subject's lifestyle, using statistical tables,
to derive the posterior probability of HIV infection.
This analysis assumes that the reference class is the set of people
with the same lifestyle.
On the other hand, if one refuses to relax the consideration
that the particular subject is either already infected or not infected,
the probability of HIV infection cannot be defined (except trivially,
that is, ``either $0$ or $1$'').

In view of these consequences of our options, we elect to include (or imply)
the clause ``disregarding the true value of $\theta$''
in the interpretation of a pivotal random variable, after the outcome is known,
so that we prevent the collapse into a trivial reference class.
A deliberate omission of this clause
would amount to voiding the \pdf\ of $E_{1\,\{x_1{=}s\}}$.

\subsubsection{The impossibility of a fiducial argument}
\label{sec:fiducial}
Fisher has emphasised repeatedly that
the fiducial assignment of a \pdf\ to the parameter
requires a carefully considered specification
of the reference class.
(He uses the terms `aggregate', `population', and `reference set',
as cited in \cite{Zabell}.)
However, when we follow this advice,
we find that the fiducial argument is not sustainable.

Suppose (for the moment) that we have established the assignment
of a \pdf\ to $E_{1\,\{x_1{=}s\}}$;
we are not thereby justified
to pair it with a corresponding \pdf\ for $\theta$,
because it would conflict with the clause
``disregarding the true value of $\theta$''
which is implied in the specification of the reference class.
A striking incongruity looms in the sentence
``the probability of $\theta$ being between $t_1$ and $t_2$,
disregarding the value of $\theta$, is $0.95$''.

Here we divorce direct predictive inference
from the fiducial argument, in terms of logical connection,
yet it can be said that the two schemes are related in intention.
In the words of A.P.~Dempster \cite{Dempster64},
``fiducial probabilities are intended
for post-data predictive interpretation''.
Also F.~Hampel \cite{Hampel}
focuses on the predictive role of fiducial probabilities.
We outline Fisher's position in Sec.~\ref{sec:discussion}.

\subsubsection{Determination of post-data probability for additive noise}
\label{sec:proof}
In relation to the issue
``what is the \pdf\ of $E_{1\,\{x_1{=}s\}}$'',
one may also ask
``how is this \pdf\ updated if, on the next day, we learn
that $\theta$ had been a random outcome of some process,
with \pdf\ $\pi(\theta)$''.
The idea is that {\em Bayesian} updating applies in this case,
as if the object of the measurement were the value of the noise, $e_1$,
and the direct information about $\theta$
were only part of the measurement process.\footnote
{{\normalsize \hspace{0.1 em}
In the words of Jaynes:
``But a telescope maker might see it differently.
For him, the errors it produces are the objects of interest to study,
and a star is only a convenient fixed object
on which to focus his instrument for the purpose of determining those errors.
Thus a given data set might serve two entirely different purposes;
one man's `noise' is another man's `signal'$\,$''
\cite{Jaynes} Ch.~7.
(Also in Ch.~8 he observes the mathematical ``reciprocity''
between random variable $\Theta$ and any ancillary random variable.)
}}
The \pdf\ of $E_{1\,\{x_1{=}s\}}$
will be identified with the prior of the alternate Bayesian treatment
(since it is meant to apply when we lack any direct information about $\theta$).

To present this argument clearly,
let us denote by ${\mathcal H}$ the assumption that
$\theta$ is the result of a random process,
corresponding to a random variable $\Theta$, of \pdf\ $\pi(\theta)$.
If ${\mathcal H}$ is accepted,
then use of Bayes' theorem is justified for probability update,
assuming an outcome $x_1$.

In our case we also have the assumption of a location measurement.
It is expressed in Eq.~\ref{eq:E_ind_Theta},
which states the requirement that
the pre-measurement \pdf\ of the noise ($E_1$)
be independent of
[any hypothetical value that might be supposed of] $\theta$.
Note that the converse property is also true, trivially:
\begin{equation} \label{eq:Theta_ind_E}
(\forall \, e) \ \ \ 
{\Pr}_{\Theta}(\theta|
{\mathcal H},\textrm{``}e_1{=}e\textrm{''})
\ \equiv \ 
{\Pr}_{\Theta}(\theta|{\mathcal H})
\ \equiv \ \pi(\theta) \, \textrm{.}
\end{equation}
That is, the prior \pdf\ associated with $\theta$
is independent of
[any hypothetical magnitude that might be supposed of]
the error of this measurement.

Note the formal symmetry
between $(\theta,{\mathcal H})$ and $(e_1,{\mathcal B}_1)$,
by comparing Eq.~\ref{eq:E_ind_Theta} with
Eq.~\ref{eq:Theta_ind_E}.\footnote
{{\normalsize \hspace{0.1 em}
This symmetry may be clouded because of certain properties
that typically are desired for $f_1(e_1)$
without being required of it:
we prefer that it average to zero, and that it also be symmetrical about zero;
moreover, it is convenient that it follow the normal distribution.
However, these properties being nonessential,
there is no real issue here.
Besides, when we want to perform a direct zero calibration of the apparatus,
we usually select a $\pi(\theta)$ having the above properties.
}}
Therefore there are two Bayesian ways of deriving
${\Pr}_{E_1}(e|{\mathcal B}_1,{\mathcal H},
\textrm{``}x_1{=}s\textrm{''})$.
As a consistency check, let us compare the results
of the two corresponding treatments.
\vspace {0.5 em} 

\noindent {\bf A}.~~ {\em The usual Bayesian treatment}
\vspace {0.5 em} \nobreak

In the usual treatment
we first update the \pdf\ of $\Theta$,
from $\pi(\theta)$ to the corresponding posterior \pdf.
We apply the familiar Bayesian formula
\[
\textrm{posterior}\ \textit{pdf}\ = \ \textrm{prior}\ \textit{pdf}
\ \times\ \textrm{likelihood}\ \times\ \textrm{normalising constant.}
\]
The likelihood function for $\theta$, given ``$x_1{=}$ \< $s$'',
is defined up to an unimportant factor:
\begin{equation} \label{eq:likel_theta_1}
{\mathrm L}_{\theta \,\{x_1{=}s\}}(t)
\ \propto \ 
{\Pr}_{X_1}(s | 
{\mathcal B}_1,\textrm{``}\theta{=}t\textrm{''}) \, \textrm{.}
\end{equation}
Considering the transformation from $X_{1\,\{\theta{=}t\}}$ to
$E_{1\,\{\theta{=}t\}}$ $\equiv$ $X_{1\,\{\theta{=}t\}}{-}$ \< $t$,
of Jacobian determinant $1$,
we obtain
\begin{equation} \label{eq:likel_theta_2}
{\Pr}_{X_1}(s|
{\mathcal B}_1,\textrm{``}\theta{=}t\textrm{''})
\ \equiv \ 
{\Pr}_{E_1}(s{-}t|
{\mathcal B}_1,\textrm{``}\theta{=}t\textrm{''}) \, \textrm{.}
\end{equation}
From Eq.s~\ref{eq:likel_theta_1}, \ref{eq:likel_theta_2},
and \ref{eq:E_ind_Theta}, we obtain
\begin{equation} \label{eq:likel_theta_3}
{\mathrm L}_{\theta \,\{x_1{=}s\}}(t)
\ \propto \ 
f_1(s{-}t) \, \textrm{.}
\end{equation}
Note that the definition of the likelihood function
does not involve ${\mathcal H}$.
(That is, with a different prior \pdf\ for $\theta$,
or in default of any prior \pdf, the likelihood function would be the same.)

The posterior \pdf\ of $\Theta$ is\footnote
{{\normalsize \hspace{0.1 em}
Strictly speaking,
if $f_1(\cdot)$ is smooth,
Bayesian updating cannot be based on the acceptance of
``$x_1{=}$ \< $s$''
because the probability of that occurence is zero
for all values of $\theta$.
Instead of an exact value for $x_1$,
we consider some small interval including that value.
In the first-order approximation
we reckon the (conditional on $\theta$) probability of this interval
as the product
${\Pr}_{X_1}(s | 
{\mathcal B}_1,$ \< ``$\theta{=}t$''$) \times$ $\Delta x_1$.
Therefore if the interval is small enough we apply
the likelihood function defined in Eq.~\ref{eq:likel_theta_1}
as the first-order approximation.
In fact it is never possible to record $x_1$ exactly;
it is always registered as a digitised entry,
which is equivalent to some interval.
}}
\[
{\Pr}_{\Theta}(t|
{\mathcal H},{\mathcal B}_1,\textrm{``}x_1{=}s\textrm{''})
\ \propto \ 
{\mathrm L}_{\theta \,\{x_1{=}s\}}(t) \; \pi(t)
\ \propto \ 
f_1(s{-}t) \; \pi(t) \, \textrm{.}
\]
Now we take advantage of the transformation from $\Theta$ to $E_1$
which is defined by $E_1 \equiv$ $s{-}$ \< $\Theta$,
of Jacobian determinant $1$, to derive
\begin{equation} \label{eq:post_e1_A}
{\Pr}_{E_1}(e|
{\mathcal B}_1,{\mathcal H},\textrm{``}x_1{=}s\textrm{''})
\ \equiv \ 
{\Pr}_{\Theta}(s {-} e|
{\mathcal B}_1,{\mathcal H},\textrm{``}x_1{=}s\textrm{''})
\ \propto \ 
f_1(e) \; \pi(s {-} e) \, \textrm{.}
\end{equation}
\vspace {0.5 em} 

\noindent {\bf B}.~~ {\em The ``instrument maker's'' Bayesian treatment}
\vspace {0.5 em} \nobreak

In the alternate Bayesian treatment
we update the \pdf\ of ${E_1}$,
from $f_1(e_1)$ to the corresponding posterior \pdf.
Again we make use of a likelihood function
but now it is the likelihood function for
${e_1}$, given ``$x_1{=}$ \< $s$''.
In analogy with the previous treatment,
we define the likelihood function
for $e_1$, given ``$x_1{=}$ \< $s$'',
without reference to ${\mathcal B}_1$:
\[
{\mathrm L}_{e_1 \,\{x_1{=}s\}}(e)
\ \propto \ 
{\Pr}_{X_1} (s |
{\mathcal H},\textrm{``}e_1{=}e\textrm{''})
\]
Considering the transformation from $X_{1\,\{e_1{=}e\}}$ to
$\Theta \equiv$ $X_{1\,\{e_1{=}e\}}{-}$ \< $e$,
we obtain
\[
{\Pr}_{X_1} (s |
{\mathcal H},\textrm{``}e_1{=}e\textrm{''})
\ \equiv \ 
{\Pr}_{\Theta}(s{-}e|
{\mathcal H},\textrm{``}e_1{=}e\textrm{''})
\ \equiv \ 
\pi(s{-}e)
\]
(in the second step we have taken into account Eq.~\ref{eq:Theta_ind_E})
so that the likelihood function is
\[
{\mathrm L}_{e_1 \,\{x_1{=}s\}}(e)
\ \propto \ 
\pi(s{-}e)
\]

Now we can derive the posterior \pdf\ of ${E_1}$,
up to a normalisation factor:
\begin{equation} \label{eq:post_e1_B}
{\Pr}_{E_1}(e|
{\mathcal B}_1,{\mathcal H},\textrm{``}x_1{=}s\textrm{''})
\ \propto \ 
{\mathrm L}_{e_1 \,\{x_1{=}s\}}(e) \; f_1(e)
\ \propto \ 
\pi(s{-}e) \; f_1(e) \, \textrm{.}
\end{equation}
\vspace {0.5 em} 

As expected, Eq.~\ref{eq:post_e1_B} is equivalent with Eq.~\ref{eq:post_e1_A},
so that the ``instrument maker's'' version of Bayesian updating
is checked as accurate.
The important point is, the prior \pdf\ in this procedure,
which is meant to apply
as long as no direct information about $\theta$ is ({\em yet}) available,
is just $f_1(e_1) \,$.\footnote
{{\normalsize \hspace{0.1 em}
If there is any doubt whether the \pdf\ for the error is legitimate
when we do not know of any \pdf\ for $\theta$,
let us refer to the symmetrical situation,
when we admit a prior \pdf\ for the parameter 
regardless of the properties of the measuring apparatus,
even regardless of whether there will be any measurement at all.
}}
In other words, we have determined the \pdf\ of $E_{1\,\{x_1{=}s\}}$,
and it turns out the same as the pre-measurement \pdf\ of $E_1$.
We have shown that $x_1$ does not specify
any recognizable subset of the reference class
that is indicated by the clause
``disregarding the unknown value of the parameter'',
with regard to random variable $E_1$.

\subsubsection{Justification of the direct pivotal argument}

A corollary of the irrelevance of $x_1$
to the \pdf\ of $E_1$
is that $x_1$ is also irrelevant to the \pdf\ of $D$
(defined in Sec.~\ref{sec:fundamental}).
This is due to the identities
\[
d \ \equiv \ x_2{-}x_1
\ \equiv \ e_2{-}e_1 \, \textrm{,}
\]
so that
\[
D \ \equiv \ E_2{-}E_1 \, \textrm{.}
\]
We have already seen that
the \pdf\ of $E_1$ is unaffected by our knowing $x_1$,
and of course so is the \pdf\ of $E_2$.
Consequently, so is the \pdf\ of $D$.
In notation,
\begin{equation}
{\Pr}_D(d|{\mathcal B}_1,{\mathcal B}_2,
\textrm{``}x_1{=}s\textrm{''})
\ \equiv \ 
{\Pr}_D(d|{\mathcal B}_1,{\mathcal B}_2)
\end{equation}

We have solved the problem stated in Sec.~\ref{sec:fundamental},
regarding the post-measurement \pdf\ of $D$.
In this way we have provided the foundation
for the direct pivotal argument,
so that we can produce the \pdf\ for
``$x_2$ given $x_1$,
disregarding the unknown value of $\theta$''.

To obtain a concrete result, note that the \pdf\ of $D$
is the marginal \pdf\ of $E_2{-}$ \< $E_1$
for any value of $\theta$,
\begin{equation}
{\Pr}_D(d|{\mathcal B}_1,{\mathcal B}_2)
\ \equiv \ 
\int \D \, \xi \; f_1(\xi) \, f_2(d{+}\xi)
\end{equation}
where $\xi$ is a dummy variable.
Considering the transformation from
$D_{\{x_1{=}s\}}$ to $X_{2\,\{x_1{=}s\}}$, we conclude that
\begin{equation}
{\Pr}_{X_2}(x_2|
{\mathcal B}_1,{\mathcal B}_2,\textrm{``}x_1{=}s\textrm{''})
\ \equiv \ 
\int \D \, \xi \; f_1(\xi) \,
f_2(x_2{-s}{+}\xi) \, \textrm{.}
\end{equation}

This result is based on disregarding the true value of the parameter.
By coincidence, fiducial prediction also results to the same \pdf,
which also coincides with the predictive \pdf\ based
on a uniform prior density for the parameter.

\section{Discussion} \label{sec:discussion}
Although the ``direct'' pivotal argument
applies only with location measurements
(a special case, even if not too uncommon)
the importance of this analysis lies in showing
an example of non-paramet\-ric predictive inference
based on parametric models.

In another paper we shall extend this result in two ways:
the predicted outcome need not be related to a location measurement,
and the prediction may be based on any number of location measurements.
However modest those developments appear in relation to the general case,
the issues raised by them require careful treatment,
so that they cannot be addressed in a short paper like the present one.

Here is a note regarding the distinction
between direct and fiducial prediction.
Fisher has not overlooked that
the problem of fiducial predictive inference based on datum $x_1$
can be solved ``directly'',
that is, not only
``after the [...] distribution of the population parameter[...]
has been obtained'' (\cite{Fisher35}, Sec.~II);
in other words:
``without discussing the possible values of the parameter $\theta$''
(\cite{Fisher56}, Sec.~V.3).
Yet he defines fiducial prediction
as derived from fiducial probability of the parameter values;
consequently the simplification he mentions is only a secondary issue.
In the present work predictive inference is defined
in the absence of any distribution for $\theta$,
therefore the possibility to also calculate it
as if from some intuitive density function of the parameter
is fortuitous, proved in the case of location measurements
but not yet guaranteed to be generally true.

\section{Conclusions}
The error (or ``noise'') of a location measurement
corresponds to a pivotal random variable.
There is an issue regarding what is the appropriate reference class
for interpreting this random variable after the outcome is known.
We show that, if we want a non-trivial and a practically useful result,
the reference class must be specified by the clause
``disregarding the unknown value of the parameter''.
In this way we also preserve correspondence
with the common usage of the term `probability'.
This clause prevents the application of the fiducial argument,
so that no \pdf\ for the parameter may be justified;
fiducial prediction is also voided by this clause.
However, the direct pivotal argument
remains valid.
It solves the problem of predictive inference, for location measurements,
without any intermediate parametric inferece.
In this way we have attained ``pure'' predictive inference;
that is, not involving any inductive component;
every step involves deduction only.

\end{document}